\UseRawInputEncoding

\documentclass[12pt]{article}

\usepackage{stmaryrd}
\usepackage{lmodern}
\usepackage{diagbox}
\usepackage{mathtools}
\usepackage{amssymb,amsmath}
\usepackage{graphicx}

\usepackage{mathrsfs}

\usepackage{amsfonts}
\usepackage{epsfig}
\usepackage{graphics}

\usepackage{lscape}
\usepackage{latexsym}

\usepackage[colorlinks=true,citecolor=black,linkcolor=black,urlcolor=blue]{hyperref}
\usepackage[usenames,dvipsnames]{color}

\usepackage{color}

\newtheorem{prop}{Proposition}

\newtheorem{lemma}{Lemma}
\newtheorem{thm}{Theorem}
\newtheorem{cor}{Corollary}

\def\Fcal{{\cal F}}
\def\Rcal{{\cal S}}
\def\eps{{\varepsilon}}

\def\proof{\noindent {\bf Proof }}
\def\qed{~~\vrule height8pt width4pt depth0pt}

\title{The genus distribution of cubic graphs and asymptotic number of rooted cubic maps with high genus } 

\author{
Zhicheng Gao\thanks{Research supported by Carleton University Research Development Fund}\\ 
School of Mathematics and Statistics\\
Carleton University\\
Ottawa, Ontario K1S5B6, Canada\\
{\small \texttt{zgao@math.carleton.ca}}
}

\begin{document}
\maketitle

\smallskip

\begin{abstract}
Let $C_{n,g}$ be the number of rooted cubic maps with $2n$ vertices on the orientable surface of genus $g$.  We show that the sequence $(C_{n,g}:g\ge 0)$ is asymptotically normal with mean and variance asymptotic to $(1/2)(n-\ln n)$ and $(1/4)\ln n$, respectively. 
We derive an asymptotic expression for $C_{n,g}$ when $(n-2g)/\ln n$ lies in any closed subinterval of $(0,2)$. 
Using rotation systems and Bender's theorem about generating functions with fast-growing coefficients, we derive simple asymptotic expressions for the numbers of rooted regular maps, disregarding the genus. In particular, we show that the number of rooted cubic maps with $2n$ vertices, disregarding the genus,   is asymptotic to $\frac{3}{\pi}\,n!6^n$.

\end{abstract}

\section{Introduction}

This paper is motivated by two lines of research. The first is about genus distributions of graphs, which is an active research area in topological graph theory 
\cite{Gro11,Gro13,GKS14,GroRobTuc89,Sta90,WanLaiLiu09}. The second is about asymptotic properties of rooted maps of high genus, which has attracted much attention recently \cite{Ang13,BudLou19,GaoKan20,Ray15,Ray21}. There is a close connection between genus distribution of graphs and enumeration of rooted maps, which we shall briefly describe below. Let $\Sigma_g$ denote the orientable surface of genus $g$. 
A {\em map} on $\Sigma_g$ is a connected graph $G$ that is {\em embedded} on $\Sigma_g$
in such a way that each component of $\Sigma_g - G$, called {\em a face}, is a simply connected region (i.e.,  a topological disk). Such an embedding is known as {\em 2-cell} (or {\em cellular}) embedding. Throughout the paper, all maps and graph embeddings are cellular embeddings on orientable surfaces.
A map on $\Sigma_g$  will be called a map with genus $g$.  A map is called {\em cubic} if all its vertices have degree 3. The dual of a cubic map is known as a {\em triangular map} (or {\em triangulation}),  whose faces all have degree 3.
Throughout the paper, a map is always {\em rooted}, meaning that a vertex and an edge incident to it are distinguished. The notion of rooted map was first introduced by Tutte  in his seminal paper \cite{Tutte-planarmap} on enumeration of planar maps. We emphasize that there is no labeling (neither vertex nor edge) in  rooted maps. Rooting a map trivializes the automorphism group which makes enumeration easier. 

On the other hand, in topological graph theory \cite{GroFur87,GT87}, one trivializes the automorphism group of an embedded graph by labeling and orienting the edges. A 2-cell embedding can be described combinatorially in terms of a permutation $\phi$ of the ends of the edges, known as {\em rotation system}.  The cycles of $\phi$ describe the cyclic ordering (say, clockwise) of the (ends of) edges at each vertex.  For each rooted map with $n$ edges, there are $(n-1)!2^{n-1}$ ways to label and to orient the $n-1$ non-root edges of the underlying graph. This gives the following connection between rooted maps and rotation systems \cite[Lemma~2.3]{Jac87}.
\begin{prop}\label{prop1}
 Each rooted map with $n$ edges corresponds to exactly  \\ 
 $ 2^{n-1}(n-1)!$   rotation systems.
\end{prop}
 For a given family $\Fcal$ of edge-labeled directed graphs, let $F_{n,g}$  denote the number of embeddings of the graphs in $\Fcal$ of $n$ edges and genus $g$. The sequence $(F_{n,g},g\ge 0)$ is called the {\em genus distribution} of $\Fcal$.  By Proposition~\ref{prop1}, $F_{n,g}2^{1-n}/(n-1)!$ is the number of rooted maps whose underlying graphs are in $\Fcal$. Genus distributions of graphs were first introduced by Gross and Furst \cite{GroFur87}.  Jackson \cite{Jac87}  derived a simple recursion for the number of rooted {\em bouquets} (maps with only one vertex, whose duals are  known as {\em unicellular maps}) using rotation system and group characters.  Using Jackson's result (also known as Harer-Zagier's formula \cite{HZ86}), Gross et al. \cite{GroRobTuc89} showed that the genus distribution of bouquets is strongly unimodal, and Stahl \cite{Sta90} derived an asymptotic relation between the number of rooted bouquets and the Stirling cycle numbers. It follows from Stahl's result (and the well-known fact that the distribution of Stirling cycle numbers is asymptotically normal) that the genus distribution of bouquets is asymptotically normal.  An exact expression for the genus distribution of bouquets  was given in \cite{Jac95}. 
 
Algorithms for computing  the genus distributions of some families of cubic  graphs  were studied in \cite{Gro11,Gro13,GKS14}. It was shown in \cite{WanLaiLiu09} that the genus distribution of a general graph can be computed using the genus distributions of some cubic graphs.

Much work has been done on the enumeration of rooted maps of constant genus since Tutte's pioneering work on planar maps \cite{Tutte-planarmap}. Many families of non-planar maps have also been enumerated; see, e.g.,  \cite{BenCan86,BenCanRic93,Gao93}.  In contrast, enumerative results on maps with high genus (both the genus and  the number of edges go to infinity) are rare.  It is clear that studying the genus distribution involves enumerating maps of high genus. Properties of large unicellular maps of high genus were  studied in \cite{Ang13,Ray15,Ray21}.  

Let $C_{n,g}$ be the number of rooted cubic maps with $2n$ vertices and genus $g$, and $J_{n,f}$ be the number of rooted cubic maps with $2n$ vertices and $f$ faces.   It follows from Euler's formula that
\begin{align}
f&=n+2-2g,\label{Euler} \\ 
C_{n,g}&=J_{n,n+2-2g}. \label{eq:face}
\end{align}
By Proposition~\ref{prop1},  $(C_{n,g}2^{3n-1}(3n-1)!: g\ge 0)$ is the the genus distribution of cubic graphs. The sequence
$(J_{n,f}2^{3n-1}(3n-1)!: f\ge 1)$ is known as the {\em region distribution}  of cubic graphs.

More recently properties of large triangulations (duals of cubic maps) of high genus were studied in \cite{BudLou19}, where an asymptotic formula, which is accurate up to a sub-exponential factor, was obtained. An asymptotic formula for such triangulations was reported in \cite{GaoKan20}; however, there is a gap in the proof. In this paper we will derive an asymptotic formula for the number of rooted cubic maps with high genus, which is accurate up to a constant factor. As in \cite{GaoKan20}, we use the Goulden-Jackson recursion for the number of rooted cubic maps; however, we shall focus on the genus polynomial $\sum_{g} C_{n,g}x^g$ in the current paper instead of the generating function $\sum_{n}C_{n,g}x^n$ in \cite{GaoKan20}.

The rest of the paper is organized as follows. In Section 2, we use rotation systems and Bender's theorem about generating functions with  fast growing coefficients to derive asymptotic expressions for the numbers of all maps and regular maps, disregarding the genus. In particular, a simple asymptotic formula is obtained for the number of rooted cubic maps, disregarding the genus. Our main results about the genus distribution of cubic graphs and asymptotic number of rooted cubic maps of high genus  are  stated at the end of this section. In section~3, we analyze the asymptotic behavior of the genus polynomials of cubic graphs using the Goulden-Jackson recursion. This is the most technical part of the paper due to the complexity of the nonlinear recursion. The proofs of our main results are then completed by using the limit theorems from \cite{GaoRic92}. Section~4 concludes our paper.

\section{The results  }
Since there are $(2n-1)!$ permutations of $2n$ elements which have exactly one cycle, by  Proposition~\ref{prop1}  the total number of rooted bouquets with $n$ edges is equal to 
\[
\frac{(2n-1)!}{(n-1)!2^{n-1}}.
\]
 Proposition~\ref{prop1} can also be used to find the asymptotic number of rooted maps, disregarding the genus, in some other families. To the best of our knowledge, such asymptotic results have not appeared  in the literature. We will use the following result about the generating functions with fast growing coefficients, which is an immediate consequence of \cite[Corollary~4]{Ben75} (See also \cite[Theorem~2]{Wri70}). We say that a sequence $(a_n)$ of positive numbers {\em grows super-exponentially} if $a_n/a_{n-1}\to \infty$ as $n\to \infty$.
 \begin{prop}\label{Bender} 
 Suppose  $(a_n)$ grows super-exponentially, $a_0=1$ and
 \begin{align}\label{eq:Ben}
 \sum_{k=1}^{n-1}a_ka_{n-k}=O(a_{n-1}).
 \end{align} Then
 \[
 \left[z^n\right]\ln \left(\sum_{n\ge 0} a_n z^n\right)\sim a_n.
 \]
 \end{prop}
 
\begin{prop}\label{prop2}
\begin{itemize}
\item[(a)]  The total number of rooted maps with $n$ edges is asymptotic to
\[
 \frac{(2n)!}{(n-1)!}2^{1-n}.
\]
\item[(b)]  For each fixed integer $d\ge 2$ and as $m\to \infty$, the number of rooted $2d$-regular maps with $m$ vertices is asymptotic to
\begin{align}\label{eq:evenr}
 \frac{(2md)!}{m!(md-1)!}(2d)^{-m}\,2^{1-md}.
\end{align}
\item[(c)]  For each fixed odd integer $r\ge 3$ and as $k\to \infty$,  the number of rooted $r$-regular maps with $2k$ vertices is asymptotic to
\begin{align}\label{eq:oddr}
 \frac{(2kr)!}{(2k)!(kr-1)!}r^{-2k}2^{1-kr}.
\end{align}
\end{itemize}
\end{prop}
\proof  For part~(a), we note \cite{JV90} that the exponential generating function of (transitive) rotation systems  is given by
\[
 \ln\left(\sum_{k\ge 0}\frac{(2k)!}{k!}z^k\right).
\]
It follows from Proposition~\ref{prop1} that the total number of rooted $n$-edged maps is equal to 
\[
\frac{1}{(n-1)!}2^{1-n}n!\left[z^n\right] \ln\left(\sum_{k\ge 0}\frac{(2k)!}{k!}z^k\right).
\]
We now verify that  $a_k:=\displaystyle \frac{(2k)!}{k!}$  grows super-exponentially and satisfies \eqref{eq:Ben}. We have
\begin{align*}
\frac{a_n}{a_{n-1}}&=2(2n-1)\to \infty.
\end{align*}
For $1\le k\le n/2$, we have
\begin{align*}
\frac{a_ka_{n-k}}{a_{k-1}a_{n-k+1}}&=\frac{(2k-1)}{(2n+1-2k)}\le 1.
\end{align*}
Consequently
\begin{align*}
\sum_{k=1}^{n-1}a_ka_{n-k}&\le a_1a_{n-1}+2\sum_{2\le k\le n/2}a_ka_{n-k}\\
&\le 2a_{n-1}+ na_2a_{n-2}\\
&=O(a_{n-1}).
\end{align*}
By Proposition~\ref{Bender}, we have
\[
\left[z^n\right] \ln\left(\sum_{k\ge 0}\frac{(2k)!}{k!}z^k\right)\sim \frac{(2n)!}{n!}.
\]
Hence the  total number of rooted $n$-edged maps is asymptotic to
\[
\frac{1}{(n-1)!}2^{1-n}n!  \frac{(2n)!}{n!},
\]
which completes the proof of part~(a).

For part~(b), we first note that each $2d$-regular graph with $m$ vertices has $md$ edges, and there are 
\[
\frac{(2md)!}{(2d)^m\,m!}
\]
permutations of $2md$ elements with exactly $m$ cycles of  length $2d$. Hence the corresponding exponential generating function of (transitive) rotation systems is given by
\[
 \ln\left(\sum_{m\ge 0}\frac{(2md)!}{(2d)^m\,m!}\frac{1}{(md)!}z^{md}\right).
\]
It follows  from Proposition~\ref{prop1} (use the substitution $w:=z^d$) that the number of rooted $2d$-regular maps with $m$ vertices is equal to
\[
  \frac{1}{(md-1)!}\,2^{1-md}(md)!\left[w^{m}\right] \ln\left(\sum_{m\ge 0}\frac{(2md)!}{(2d)^m\,m!}\frac{1}{(md)!}w^{m}\right).
\]
We now verify that  $a_m:=\displaystyle \frac{(2md)!}{(2d)^m\,m!}\frac{1}{(md)!}$  grows super-exponentially and satisfies \eqref{eq:Ben}.  For each fixed $d\ge 2$ and as $m\to \infty$, we have
\begin{align*}
\frac{a_m}{a_{m-1}}&=\frac{(2md)(2md-1)\cdots (2md-2d+1)}{(2d)(m)(md)(md-1)\cdots (md-d+1)} \\
&=\frac{1}{md}2^{d-1}\prod_{1\le j\le 2d-1,2\nmid j}(2md-j)\\
&\sim 2^{2d-1}(md)^{d-1},
\end{align*}
which shows that $a_m$ grows super-exponentially.
For $1\le k\le m/2$, we  have 
\[
k\le m-k+1,~~
\frac{2kd-j}{2(m-k+1)d-j}\le 1,
\]
and hence
\begin{align*}
\frac{a_ka_{m-k}}{a_{k-1}a_{m-k+1}}&=\frac{m-k+1}{k}\prod_{1\le j\le 2d-1,2\nmid j}\frac{2kd-j}{2(m-k+1)d-j}\\
&\le \frac{m-k+1}{k}\frac{2kd-1}{2(m-k+1)d-1}\\
&=\frac{2k(m-k+1)d-(m-k+1)}{2k(m-k+1)d-k}\\
&\le 1.
\end{align*}
Now  the asymptotic expression \eqref{eq:evenr} follows from the same argument as that for part~(a).

For part~(c),  since $r\ge 3$ is odd,  each $r$-regular graph must have an even number of vertices, and each $r$-regular graph with $2k$ vertices has exactly $kr$ edges. Since there are 
\[
\frac{(kr)!}{r^{2k}\,(2k)!}
\]
permutations of $2kr$ elements with exactly $2k$ cycles of  length $r$,  the corresponding exponential generating function of (transitive) rotation systems is given by
\[
 \ln\left(\sum_{k\ge 0}\frac{(2kr)!}{r^{2k}\,(2k)!}\frac{1}{(kr)!}z^{kr}\right).
\]
It follows  from Proposition~\ref{prop1} (use the substitution $w:=z^r$) that the number of rooted $r$-regular maps with $2n$ vertices is equal to
\[
  \frac{1}{(nr-1)!}\,2^{1-nr}(nr)!\left[w^{n}\right] \ln\left(\sum_{n\ge 0}\frac{(2nr)!}{r^{2n}\,(2n)!}\frac{1}{(nr)!}w^{n}\right).
\]
We now verify that  $a_k:=\displaystyle \frac{(2kr)!}{r^{2k}\,(2k)!}\frac{1}{(kr)!}$  grows super-exponentially and satisfies \eqref{eq:Ben}.  For each fixed $r\ge 3$ and as $n\to \infty$, we have
\begin{align*}
\frac{a_n}{a_{n-1}}
&=\frac{1}{r^2n(2n-1)}2^{r-1}\prod_{1\le j\le 2r-1,2\nmid j}(2nr-j)\\
&\sim 2^{2r-2}(nr)^{r-2}.
\end{align*}
For $1\le k\le n/2$, we  have 
\[
\frac{2kr-j}{2(n-k+1)r-j}\le 1,
\]
and hence (taking the terms corresponding to $j=1$ and $j=2r-1$ from the product)
\begin{align*}
\frac{a_ka_{n-k}}{a_{k-1}a_{n-k+1}}&=\frac{(n-k+1)(2n-2k+1)}{k(2k-1)}\prod_{1\le j\le 2r-1,2\nmid j}\frac{2kr-j}{2(n-k+1)r-j}\\
&\le \frac{(n-k+1)(2n-2k+1)}{k(2k-1)}\frac{2kr-1}{2(n-k+1)r-1}\frac{2(k-1)r+1}{2(n-k)r+1}\\
&=\frac{(2kr-1)(n-k+1)}{k(2(n-k+1)r-1)}\frac{(2(k-1)r+1)(2n-2k+1)}{(2k-1)(2(n-k)r+1)}\\
&=\frac{2kr(n-k+1)-(n-k+1)}{2k(n-k+1)r-k}\frac{(2(k-1)r+1)(2n-2k+1)}{(2k-1)(2(n-k)r+1)}\\
&\le \frac{2(2k-1)(n-k)r+(2k-1)-2(n+1-2k)(r-1)}{2(2k-1)(n-k)r+(2k-1)} \\
&\le 1.
\end{align*}
Now  the asymptotic expression \eqref{eq:oddr} follows from the same argument as that for part~(a).
\qed

Let $C_n=\sum_{g\ge 0}C_{n,g}$ be the total number of rooted cubic maps with $2n$ vertices. The next result follows immediately from  Proposition~\ref{prop2}(c) and Stirling's formula.
\begin{cor}\label{cor1} 
The total number of rooted cubic maps with $2n$ vertices is
\begin{align*}
C_{n}&\sim \frac{3}{\pi}\,n!\,6^{n}\sim \frac{6}{\sqrt{2\pi }}\left(\frac{6}{e}\right)^{n}n^{n+\frac{1}{2}}.
\end{align*}
\end{cor}
We would like to point out that the above asymptotic expression can also be obtained using a bijection between cubic maps and certain family of $\lambda$-terms \cite[Theorems~3.3 and 3.4]{BGJ13}.

 Our main results are summarized in the following two theorems.
\begin{thm}\label{thm2} 
\begin{itemize}
\item[(a)] The genus distribution of cubic graphs is asymptotically normal  with mean and variance, respectively, asymptotic to
$\frac{n-\ln n}{2}$ and $\frac{\ln n}{4}$.
That is,
\begin{align*}
\sum_{g\le \frac{n-\ln n}{2}+\frac{t\sqrt{\ln n}}{2} }\frac{C_{n,g}}{C_{n}}\sim \frac{1}{\sqrt{2\pi}}\int_{-\infty}^t \exp\left(-\frac{x^2}{2} \right)dx.
\end{align*} 
\item[(b)] The region distribution of cubic graphs is asymptotically normal with mean and variance both asymptotic to
 $\ln n$.
\end{itemize}
\end{thm}

\begin{thm}\label{thm3} 
Let $\eps$ be any small positive constant. There is a function $K(y)$ which  is analytic in $[\eps,2-\eps]$ such that
\begin{align} \label{eq:Casy}
C_{n,g}&\sim \frac{\sqrt{2}}{3}K\left(\frac{n-2g}{\ln n} \right)\left(\frac{\ln n}{n-2g}\right)^2\, 6^n\frac{(n-1)!}{(n-2g)!} (\ln n)^{n-2g},
\end{align}
uniformly for all  $g,n\to \infty$ satisfying
\[
\frac{n-2g}{\ln n}\in [\eps,2-\eps].
\]
\end{thm}

\section{Proofs of Theorems~1 and 2}

Define
\begin{align}\label{eq:HC}
H_{n,g} &= (3n+2)C_{n,g} \quad\mbox{for}\quad n\ge1, 
\end{align}
\begin{align}\label{eq:Hinitial}
H_{-1,g}=\llbracket g=0\rrbracket/2,\quad H_{0,g}=2\llbracket g=0\rrbracket,\quad H_{n,-1}=0.
\end{align}
Goulden and Jackson \cite{GJ08} derived the following recursion for $(n,g)\ne (-1,0)$:
\begin{align} 
H_{n,g} &= \frac{4n(3n+2)(3n-2)}{n+1}H_{n-2,g-1} \nonumber \\
&~~+  \frac{4(3n+2)}{n+1}\sum_{k=-1}^{n-1}\sum_{h=0}^gH_{k,h}H_{n-2-k,g-h}.\label{eq:Hgn}
\end{align}

Define
\begin{align*}
H_n(x)&=\frac{3n+2}{n!}6^{-n}\sum_{g\ge 0}C_{n,g}x^g=\frac{1}{n!}6^{-n}\sum_{g\ge 0}H_{n,g}x^g,\\
J_n(y)&=\frac{3n+2}{n!}6^{-n}\sum_{f\ge 1}J_{n,f}y^{f}.
\end{align*}
Using \eqref{eq:face}, we obtain 
\begin{align}
J_n(y)&=\frac{3n+2}{n!}6^{-n}\sum_{g\ge 0}C_{n,g}y^{n+2-2g} \nonumber  \\
&=H_n(1/y^2)y^{n+2}. \label{eq:JH}
\end{align}
It follows from \eqref{eq:Hinitial}  and \eqref{eq:Hgn}  that
\begin{align}
H_0(x)&=2, \label{eq:Hx0}\\
H_1(x)&=\frac{20+5x}{6}, \label{eq:Hx1}\\
H_2(x)&=\frac{32+28x}{9}, \label{eq:Hx2}\\
H_n(x)&=\frac{2(3n+2)}{3n(n+1)}H_{n-1}(x) +\frac{9n^2-4}{9(n^2-1)}xH_{n-2}(x)&\nonumber \\
&~~+\frac{3n+2}{9n(n^2-1)}\sum_{k=0}^{n-2}\frac{1}{{n-2\choose k}}H_{k}(x)H_{n-2-k}(x).&(n\ge 2)\label{eq:Hxn}
\end{align}
Setting 
$
x=1/y^2
$
and using \eqref{eq:JH}--\eqref{eq:Hxn}, we obtain
\begin{align}
J_1(y)&=\frac{y(20y^2+5)}{6}, \label{eq:J1}\\
J_2(y)&=\frac{4y^2(8y^2+7)}{9}, \label{eq:J2}\\
J_3(y)&= \frac{11y(336y^4 + 664y^2 + 105)}{1296}, \label{eq:J3}\\
J_4(y)&= y^2\left(\frac{448}{243}y^4+\frac{1631}{243}y^2+\frac{1183}{324}\right), \label{eq:J4}\\
J_5(y)&= \frac{17y(27456y^6 + 163248y^4 + 198396y^2 + 25025)}{466560}, \label{eq:J5}\\
J_6(y)&=\frac{y^2\left(3072y^6 + 27532y^4 + 61185y^2 + 26261\right)}{6561}, \label{eq:J6}
\end{align}
and for $n\ge 7$ (separating the six terms corresponding to $k\in\{0,1,2,n-2,n-3,n-4\}$ in the last summation of \eqref{eq:Hxn} )
\begin{align}
J_n(y)&=\frac{2(3n+2)y}{3n(n+1)}J_{n-1}(y) \nonumber \\
&~~+\left(\frac{9n^2-4}{9(n^2-1)}+\frac{4(3n+2)y^2}{9n(n^2-1)}\right)J_{n-2}(y) \nonumber\\
&~~+\frac{2(3n+2)}{9n(n^2-1)(n-2)}J_1(y)J_{n-3}(y) \label{eq:JRec}\\
&~~+\frac{4(3n+2)}{9n(n^2-1)(n-2)(n-1)}J_2(y)J_{n-4}(y)\nonumber\\
&~~+\frac{3n+2}{9n(n^2-1)}\sum_{k=3}^{n-5}\frac{1}{{n-2\choose k}}J_k(y)J_{n-2-k}(y). \nonumber
\end{align}

To derive asymptotics of $J_n(y)$, we first show  that $n^{-y}J_n(y)=O(1)$ when $0\le y\le 2$.
\begin{lemma}  \label{Lem1}
Let
\begin{align}\label{eq:hJ0}
h_n(y)&:=n^{-y}J_n(y).
\end{align}
For $n\ge 2$ and $0\le y\le 2$, we have
\begin{align}
h_n(y)&\le 9n, &  \label{eq:h1}\\
h_n(y)&\le \exp\left(10-10/n\right).  \label{eq:hSmall} 
\end{align}
\end{lemma} 
{\bf Proof}   Recursion \eqref{eq:JRec} can be rewritten as
\begin{align}
h_n(y)&=\frac{2(3n+2)y}{3n(n+1)}\left(\frac{n-1}{n}\right)^{y}h_{n-1}(y) \nonumber \\
&~~+\left(\frac{9n^2-4}{9(n^2-1)}+\frac{4(3n+2)y^2}{9n(n^2-1)}\right)\left(\frac{n-2}{n}\right)^{y}h_{n-2}(y) \nonumber\\
&~~+\frac{2(3n+2)}{9n(n^2-1)(n-2)}\left(\frac{n-3}{n}\right)^{y}J_1(y)h_{n-3}(y) \label{eq:hRec}\\
&~~+\frac{4(3n+2)}{9n(n^2-1)(n-2)(n-3)}\left(\frac{n-4}{n}\right)^{y}J_2(y)h_{n-4}(y)\nonumber\\
&~~+\frac{3n+2}{9n(n^2-1)}\sum_{k=3}^{n-5}\frac{1}{{n-2\choose k}}\left(\frac{k(n-2-k)}{n}\right)^{y}h_k(y)h_{n-2-k}(y). \nonumber
\end{align}

Using \eqref{eq:J1}--\eqref{eq:J6} and \eqref{eq:hRec}, it is easy to check (using Maple) that both \eqref{eq:h1} and \eqref{eq:hSmall} hold for $2\le n\le 19$. 
So we assume $n\ge 20$ and move on to the inductive step.

Using
\begin{align*}
\frac{3n+2}{3(n+1)}&\le 1,\\
\frac{9n^2-4}{9(n^2-1)}+\frac{4(3n+2)y^2}{9n(n^2-1)}&\le 1+\frac{5}{9(n^2-1)}+\frac{4y^2}{3n(n-1)},\\
\end{align*}
and dividing both sides of \eqref{eq:hRec} by $9n$, we obtain 
\begin{align}
\frac{h_n(y)}{9n}&\le 
\frac{2y}{n}\left(\frac{n-1}{n}\right)^{y+1}+\left(\frac{n-2}{n}\right)^{y+1}+\frac{5}{9(n^2-1)} +\frac{4y^2}{3n(n-1)}\nonumber\\
&~~+\frac{2}{3n(n-1)}\left(\frac{J_1(y)}{n-2}+\frac{2J_2(y)}{(n-2)(n-3)}\right)\nonumber\\
&~~+\frac{3}{n(n-1)}\sum_{k=3}^{n-5}\frac{1}{{n-2\choose k}}\left(\frac{k(n-2-k)}{n}\right)^{y+1}. \label{eq:hsnTerm4}
\end{align}

Using $n\ge 20$, $0\le y\le 2$,
\begin{align}
(1-t)^s&\le e^{-st}\le  1-st+\frac{s^2t^2}{2},\quad\quad (0<t<1, s\ge 0) \label{eq:Ineq}\\
\frac{J_1(y)}{n-2}+\frac{2J_2(y)}{(n-2)(n-3)}&=\frac{y(20y^2+5)}{6(n-2)}+\frac{8y^2(8y^2+7)}{9(n-2)(n-3)}<\frac{5}{2}, \label{eq:J12bound}
\end{align}
 we obtain
\begin{align}
\frac{h_n(y)}{9n}
&\le \frac{2y}{n}+\left(1-\frac{2(y+1)}{n}+\frac{2(y+1)^2}{n^2}\right) +\frac{5}{9(n^2-1)}+\frac{16+5}{3n(n-1)}\nonumber \\
&~~+\frac{3}{n(n-1)}\sum_{k=3}^{n-5}\frac{1}{{n-2\choose k}}\left(\frac{k(n-2-k)}{n}\right)^{y+1}. \label{eq:check1}
\end{align}

With the help of Maple, it can be checked that, for $n\ge 20$ and $0\le y\le 2$,  
\begin{align*}
\sum_{k=3}^{n-5}\frac{1}{{n-2\choose k}}\left(\frac{k(n-2-k)}{n}\right)^{y+1}
&< \frac{1}{4}.
\end{align*}
Substituting this into \eqref{eq:check1}, we obtain 
\begin{align*}
\frac{h_n(y)}{9n}
&\le 1-\frac{2}{n}+\frac{1}{19n}\left( 18+\frac{5}{9}+\frac{21}{3}+\frac{3}{4}\right)\\
&<1.
\end{align*}
This completes the proof of \eqref{eq:h1}.

The proof of \eqref{eq:hSmall} is similar.  Dividing both sides of \eqref{eq:hRec} by $\exp(10-10/n)$, 
 and replacing $h_k(y)$ in the last line of  \eqref{eq:hRec} by $9k$,  we obtain
\begin{align}
\frac{h_n(y)}{\exp(10-10/n)}
&\le \frac{2y}{n}\left(\frac{n-1}{n}\right)^y+\left(\frac{n-2}{n}\right)^y\exp\left(-\frac{20}{n(n-2)}\right)  \nonumber \\
&~~+\frac{5}{9(n^2-1)}+\frac{4y^2}{3n(n-1)} \nonumber \\
&~~+\frac{2}{3n(n-1)} \left(\frac{J_1(y)}{n-2}+\frac{2J_2(y)}{(n-2)(n-3)}\right) \label{yRec}\\
&~~+\frac{3}{n(n-1)}\sum_{k=3}^{n-5}\frac{k}{{n-2\choose k}}\left(\frac{k(n-2-k)}{n}\right)^{y}. \nonumber
\end{align}
Using  \eqref{eq:Ineq},  \eqref{eq:J12bound},  and
\begin{align*}
\exp\left(-\frac{20}{n(n-2)}\right)&\le 1-\frac{4}{5}\frac{20}{n(n-2)}, &(n\ge 20)\\ 
\sum_{k=3}^{n-5}\frac{k}{{n-2\choose k}}\left(\frac{k(n-2-k)}{n}\right)^{y}
&< \frac{1}{4},  &(n\ge 20,0\le y\le 2.)
\end{align*}
we obtain from \eqref{yRec} that
\begin{align*}
\frac{h_n(y)}{\exp(10-10/n)}
&\le \frac{2y}{n}\left(1-\frac{y}{n}+\frac{y^2}{2n^2}\right)+\left(1-\frac{2y}{n}+\frac{2y^2}{n^2}\right)\left(1-\frac{16}{n(n-2)}\right)  \\
&~~+\frac{5}{9(n^2-1)}+\frac{21}{3n(n-1)} +\frac{3}{4n(n-1)}\\
&\le  1-\frac{1}{n(n-2)} \left(16-7-\frac{8}{20}-\frac{64}{20}-\frac{5}{20}-\frac{3}{4}\right)\\
&\le 1. 
\end{align*}
 This completes the proof of \eqref{eq:hSmall}. \qed

\medskip

\noindent {\bf Remark} 
The range of $y$ in Lemma~\ref{Lem1} can be extended. In fact we believe that $h_n(y)$ is bounded by a function of $y$ for all $y\ge 0$, but we are unable to prove this at this stage.

\medskip

The next result gives asymptotics for $h_n(y)$ which will be used to prove our main results (Theorems~1 and 2). 
In what follows, $\Re(z)$ and $\Im(z)$ stand for the real and imaginary parts of a complex number $z$, respectively, and all the big-$O$ terms are independent of $y$. 
Define
\[
\Rcal:=\{y\in {\mathbb C}:|y|\le 2, \Re(y)\ge 0,|y|-\Re(y)<1\}.
\] 
We shall also use the Iverson bracket $\llbracket P \rrbracket $, which is equal to 1 if the predicate $P$ is true and 0 otherwise.
\begin{thm} \label{thm4} 
There is a function $K(y)$ such that  in $\Rcal$  $K(y)$  is analytic and
\begin{align}
h_n(y)&= K(y)+O\left(n^{-(1+\Re(y)-|y|)}+n^{-2\Re(y)}\ln n\right) +\llbracket \Re(y)=0\rrbracket O(1). \label{eq:hsAsy}
\end{align}
\end{thm}
\proof  
We have
\begin{align}
|h_n(y)|&=\left|  J_n(y)\right| n^{-\Re(y)}\le J_n(|y|) n^{-\Re(y)}= h_n(|y|)n^{|y|-\Re(y)}. \label{eq:abs}
\end{align}
Using Lemmas~1 and 2, we have
\begin{align}
|h_n(y)|&\le e^{10}n^{|y|-\Re(y)}. &(|y|\le 2) \label{eq:hBoundUniform}
\end{align}
It follows from  \eqref{eq:hRec} that, for $|y|\le 2$,
\begin{align}
h_n(y)-h_{n-2}(y)&=\frac{2y}{n}\left(h_{n-1}(y)-h_{n-2}(y)\right)+O\left(n^{-(2+\Re(y)-|y|)}\right),\label{eq:hRec1}
\end{align}
or 
\begin{align}
h_n(y)-h_{n-1}(y)&=\left(\frac{2y}{n}-1\right)\left(h_{n-1}(y)-h_{n-2}(y)\right)+O\left(n^{-(2+\Re(y)-|y|)}\right).\label{eq:hRec2}
\end{align}
Define
\begin{align*}
d_n(y):=|h_n(y)-h_{n-1}(y)|.
\end{align*}
Using \eqref{eq:hRec2} and 
\begin{align*}
\left|\frac{2y}{n}-1\right|
&=\left(\left(\frac{2\Re(y)}{n}-1\right)^2+\left(\frac{2\Im(y)}{n}\right)^2\right)^{1/2}\\
&=1-\frac{2\Re(y)}{n}+O\left(\frac{1}{n^2}\right),
\end{align*}
we obtain 
\begin{align*}
d_n(y)&\le \left(1-\frac{2\Re(y)}{n}\right)d_{n-1}(y)+O\left(n^{-(2+\Re(y)-|y|)}\right).
\end{align*}
It is known (use  ``rsolve'' in Maple or \cite[Section 2.2]{Mic90}) that the solutions to the recursion 
\[
d_n= \left(1-\frac{a}{n}\right)d_{n-1}+n^{-b}\quad\quad (a\ge 0,n\ge 1)
\]
satisfy
\begin{align*}
d_n&=O\left(n^{-a}\left(1+\llbracket b-a=1\rrbracket\ln n\right)+n^{-(b-1)}\right).
\end{align*}
Hence
\begin{align*}
d_n(y)&=O\left(n^{-2\Re(y)}\left(1+\llbracket \Re(y)+|y|=1\rrbracket\ln n\right)+n^{-(1+\Re(y)-|y|)}\right).
\end{align*}
That is,
\begin{align}
h_n(y)-h_{n-1}(y)&=O\left(n^{-(1+\Re(y)-|y|)}+n^{-2\Re(y)}\ln n\right) +\llbracket \Re(y)=0\rrbracket O(1). \label{eq:dBound}
\end{align}

On the other hand,   using \eqref{eq:hRec1}  and
\[
\frac{h_{n-2}(y)}{n-1}-\frac{h_{n-2}(y)}{n}=\frac{h_{n-2}(y)}{n(n-1)}=O\left(n^{-(2+\Re(y)-|y|)}\right),
\]
we obtain
\begin{align}
&~h_n(y)-h_{n-2}(y)-2y\left(\frac{h_{n-1}(y)}{n}-\frac{h_{n-2}(y)}{n-1}\right)\nonumber\\
& =h_n(y)-h_{n-2}(y)-\frac{2y}{n}\left(h_{n-1}(y)-h_{n-2}(y)\right)+2y\left(\frac{1}{n-1}-\frac{1}{n}\right)h_{n-2}(y)\nonumber\\
&=O\left(n^{-(2+\Re(y)-|y|)}\right).\label{eq:hBound}
\end{align}
Set
\[
R_n(y):=h_n(y)-h_{n-2}(y)-2y\left(\frac{h_{n-1}(y)}{n}-\frac{h_{n-2}(y)}{n-1}\right).
\]
Since $J_n(y)$ is a polynomial, by \eqref{eq:hJ0}  $R_n(y)$ is analytic everywhere for each $n$. It follows from \eqref{eq:hBound}, Weierstrass M-test and Morera's theorem \cite{Rud} that the series
$\sum_{k\ge 3}R_n(y)=K_1(y)$ is analytic when $\Re(y)-|y|> -1$, and 
\begin{align}\label{eq:K1}
\sum_{k=3}^nR_n(y)=K_1(y)+O\left(n^{-(1+\Re(y)-|y|)}\right).
\end{align}

Summing both sides of \eqref{eq:hBound} from 3 to $n$ (noting the cancellations from the telescoping sum), we obtain
\begin{align*}
h_{n}(y)+h_{n-1}(y)-h_1(y)-h_2(y)+yh_1(y)&=K_1(y)+O\left(n^{-(1+\Re(y)-|y|)}\right),
\end{align*}
or
\begin{align}\label{eq:h12}
h_{n}(y)+h_{n-1}(y)&=K_1(y)+(1-y)h_1(y)+h_2(y)+O\left(n^{-(1+\Re(y)-|y|)}\right).
\end{align}
 Combining this with \eqref{eq:dBound},  and setting
\begin{align}\label{eq:K}
 K(y):=\frac{K_1(y)+(1-y)h_1(y)+h_2(y)}{2},
\end{align}
 we obtain \eqref{eq:hsAsy}.
\qed

\medskip

\noindent {\bf Remark} Using Theorem~\ref{thm4},  \eqref{eq:JH} and \eqref{eq:hJ0}, we have
\[
H_n(1)\sim K(1)n.
\]
Noting
\[
C_n=\sum_{g\ge 0}C_{n,g}=\frac{1}{3n+2}\sum_{g\ge 0}H_{n,g}=\frac{n!}{3n+2}6^nH_n(1),
\]
 we obtain
\begin{align*}
C_n&\sim \frac{K(1)}{3}n!6^n.
\end{align*}
Comparing this with the asymptotic expression in Corollary~1, we obtain
\[
K(1)=\frac{9}{\pi}.
\]
\medskip

\noindent {\bf Proof of Theorem~1}  
Using \eqref{eq:hJ0} and Theorem~\ref{thm4}, we have, uniformly for $x$ in a small neighborhood of 1, that 
\begin{align}
H_n(x)&\sim K(1/\sqrt{x})n^{1/\sqrt{x}}x^{(n+2)/2}.\label{eq:Jasy}
\end{align}

We would like to point out that Hwang's quasi-power theorem \cite{Hwang98} does not apply directly here because of the the factor $x^{(n+2)/2}$ appearing in  \eqref{eq:Jasy}.
 It is possible to apply the quasi-power theorem to $J_n(y)$, and then use \eqref{eq:face} to obtain (a).  This would also give the convergence rate.  However, we shall apply
\cite[Theorem~2]{GaoRic92} directly here.   In terms of the notations in \cite{GaoRic92}, we have
\begin{align*}
s&=\ln x,\\
m_{r}(s)&=\frac{d}{ds}\frac{s}{2}=\frac{1}{2},\\
m_{\alpha}(s)&=\frac{d}{ds} e^{-s/2}=-\frac{1}{2}e^{-s/2},\\
B_r(s)&=0,\\
B_{\alpha}(s)&=\frac{d}{ds}m_{\alpha}(s) =\frac{1}{4}e^{-s/2}.
\end{align*}
Now part~(a) follows from Theorem~2~(case~(2)) of \cite{GaoRic92} by noting 
\begin{align*}
m_r(0)n+m_{\alpha}(0)\ln n=\frac{1}{2}(n-\ln n),\\
B_{\alpha}(0)\ln n=\frac{1}{4}\ln n.
\end{align*}

Part~(b) follows immediately from part~(a) and \eqref{eq:face}. \qed

\medskip

\noindent {\bf Proof of Theorem~2}  We may use Theorem~\ref{thm4} and apply the standard saddle-point method directly. 
In what follows,  we shall apply \cite[Theorem~4]{GaoRic92} to $J_n(y)$. More precisely, to avoid the parity issue, we shall consider
\begin{align}\label{eq:GJ}
G_n(y)=\left(\sqrt{y}\right)^{\llbracket 2\nmid n\rrbracket}J_{n}(\sqrt{y}) .
\end{align}
By Lemma~\ref{Lem1}, $G_n(y)$ is a polynomial in $y$. For $y\ne 0$,  write
\begin{align*}
y&=\rho e^{i\theta},  &(-\pi<\theta \le \pi )\\
\sqrt{y}&=\sqrt{\rho}\exp\left(i\theta/2\right).
\end{align*}
Then
\[
\Re(\sqrt{y})=\sqrt{\rho}\cos(\theta/2)\ge 0.
\]

Using Theorem~\ref{thm4}, we obtain, for some small positive constants $\eps$ and $\delta$,
\begin{align*}
G_n(y)&\sim K(\sqrt{y})\left(\sqrt{y}\right)^{\llbracket 2\nmid n\rrbracket}n^{\sqrt{y}}, &  (\eps\le \rho\le 4-\eps, |\theta|\le \pi-\eps) \\
\frac{G_n(y)}{G_n(\rho)}&= O\left(n^{-\delta}\right). & ( \eps\le \rho \le 4-\eps,~|\theta|\ge \pi-\eps) 
\end{align*}
Applying  \cite[Theorem~4(2)]{GaoRic92} with $\alpha(s)=e^{s/2}$, we obtain
\begin{align}
[y^k]G_{n}(y)&\sim \frac{K(\sqrt{\rho})}{\sqrt{2\pi(\ln n)\sqrt{\rho}/4}}(\sqrt{\rho})^{\llbracket 2\nmid n\rrbracket}n^{\sqrt{\rho}}\rho^{- k}, \label{eq:GnkAsy}
\end{align}
where $\rho$ satisfies
\begin{align}
\frac{\sqrt{\rho}}{2}&=\frac{k}{\ln n}. \quad \quad (\eps\le \rho \le 4-\eps)\label{eq:rho} 
\end{align}
Substituting \eqref{eq:rho} into \eqref{eq:GnkAsy}, we obtain
\begin{align}
[y^k]G_{n}(y)&\sim \frac{K(2k/\ln n)}{\sqrt{k\pi}}e^{2k}
\left(\frac{2k}{\ln n}\right)^{\llbracket 2\nmid n\rrbracket -2 k}.
\end{align}
Using  \eqref{eq:face} and \eqref{eq:GJ}, we obtain
\begin{align*}
[y^k]G_{n}(y)&=\llbracket 2\mid n\rrbracket [y^{2k}]J_n(\sqrt{y}) +\llbracket 2\nmid n\rrbracket [y^{2k-1}]J_n(\sqrt{y}) \\
&=\frac{3n+2}{n!}6^{-n}\left(\llbracket 2\mid n\rrbracket C_{n,\frac{n}{2}+1-k}+\llbracket 2\nmid n\rrbracket C_{n,\frac{n+1}{2}+1-k}\right)\\
&=\frac{3n+2}{n!}6^{-n}C_{n,\frac{n+\llbracket 2\nmid n\rrbracket}{2}+1-k}.
\end{align*}
Combining this with \eqref{eq:GnkAsy} and setting $k=\frac{n+\llbracket 2\nmid n\rrbracket}{2}+1-g$, we obtain
\begin{align*}
C_{n,g}&\sim \frac{(n-1)!}{3}6^n \frac{K(2k/\ln n)}{\sqrt{k\pi}}e^{2k}
\left(\frac{2k}{\ln n}\right)^{\llbracket 2\nmid n\rrbracket -2 k}\\
&\sim  \frac{(n-1)!}{3}6^n \frac{\sqrt{2}\,K((n-2g)/\ln n)}{\sqrt{(n-2g)\pi}}e^{n-2g+2+\llbracket 2\nmid n\rrbracket}\\
&~~\quad\quad \times \left(\frac{n-2g+2+\llbracket 2\nmid n\rrbracket}{\ln n}\right)^{2g-n-2}\\
&\sim \frac{(n-1)!}{3}6^n \frac{\sqrt{2}\,K((n-2g)/\ln n)}{\sqrt{(n-2g)\pi}}
\left(\frac{n-2g}{\ln n}\right)^{2g-n-2}\\
&~~\quad\quad \times  e^{n-2g+2+\llbracket 2\nmid n\rrbracket}
\left(1+\frac{2+\llbracket 2\nmid n\rrbracket}{n-2g}\right)^{-(n-2g)}\\
&\sim \frac{(n-1)!}{3}6^n \frac{\sqrt{2}\,K((n-2g)/\ln n)}{\sqrt{(n-2g)\pi}}
\left(\frac{n-2g}{\ln n}\right)^{2g-n-2} e^{n-2g}.
\end{align*}
Now Theorem~\ref{thm3} follows by using  Stirling's formula:
\begin{align*}
(n-2g)! &\sim \sqrt{2\pi(n-2g)} \left(\frac{n-2g}{e}\right)^{n-2g} &(n-2g\to \infty).
\end{align*}
\qed

\section{Conclusion}

Using the Goulden-Jackson recursion for the number of rooted cubic maps,  we derived an asymptotic formula for the number of rooted cubic maps with $2n$ vertices and genus $g$ when $n,g\to \infty$ and  $(n-2g)/\ln n$ lies in any closed subinterval of $(0,2)$. The asymptotic formula is accurate up to a constant factor. We also showed that the genus distribution of cubic graphs is asymptotically normal with mean and variance, respectively, asymptotic to $\frac{1}{2}(n-\ln n)$ and $\frac{1}{4}\ln n$.  Asymptotic formulas were also obtained for the number of rooted regular maps, disregarding the genus, for constant degree and as the number of vertices going to infinity. 

\subsection*{Acknowledgement}
The author would like to thank Mihyun Kang for helpful comments on an earlier draft which improves the presentation. The author is also grateful to an anonymous reviewer who read the paper carefully and gave constructive suggestions.

\end{document}